\documentclass[11pt,letter]{article}

\usepackage{amsmath} 
\usepackage{etoolbox}
\usepackage{amssymb} 
\usepackage{graphicx} 
\usepackage{hyperref}
\usepackage{enumerate}
\hypersetup{bookmarks=true}
\usepackage{amscd}
\usepackage{amsthm}

\usepackage{calc}
\usepackage{tikz}
\usetikzlibrary{decorations.markings}
  
\tikzstyle{vertex}=[circle, draw, inner sep=0pt, minimum size=6pt]
\newcommand{\vertex}{\node[vertex]}

\def\qed{\hfill \hbox{\vrule width4pt depth2pt height6pt}}

\theoremstyle{plain}
\newtheorem{Lemma}{Lemma}[section]
\newtheorem{Theorem}[Lemma]{Theorem}
\newtheorem{Proposition}[Lemma]{Proposition}
\newtheorem{Corollary}[Lemma]{Corollary}

\theoremstyle{definition}

\newcommand{\Aut}{\mathop{\mathrm{Aut}}\nolimits}

\newcommand{\Sym}{\mathop{\mathrm{Sym}}\nolimits}
\newcommand{\Inn}{\mathop{\mathrm{Inn}}\nolimits}
\newcommand{\Cay}{\mathop{\mathrm{Cay}}\nolimits}

\newcommand{\supp}{\mathop{\mathrm{supp}}\nolimits}

\begin{document}

\title{Automorphism group of the complete transposition graph}

\author{Ashwin~Ganesan%
  \thanks{53 Deonar House, Deonar Village Road, Mumbai - 88, India. Correspondence address: 
\texttt{ashwin.ganesan@gmail.com}.}
}

\date{}

\maketitle


\begin{abstract}
\noindent  The complete transposition graph is defined to be the graph whose vertices are the elements of the symmetric group $S_n$, and two vertices $\alpha$ and $\beta$ are adjacent in this graph iff there is some transposition $(i,j)$ such that $\alpha=(i,j) \beta$.  Thus, the complete transposition graph is the Cayley graph $\Cay(S_n,S)$ of the symmetric group generated by the set $S$ of all transpositions.  An open problem in the literature is to determine which Cayley graphs are normal. It was shown recently that the Cayley graph generated by 4 cyclically adjacent transpositions is not normal.  In the present paper, it is proved that the complete transposition graph is not a normal Cayley graph, for all $n \ge 3$.  Furthermore, the automorphism group of the complete transposition graph is shown to equal
\[ \Aut(\Cay(S_n,S)) = (R(S_n) \rtimes \Inn(S_n)) \rtimes \mathbb{Z}_2, \]
where $R(S_n)$ is the right regular representation of $S_n$, $\Inn(S_n)$ is the group of inner automorphisms of $S_n$, and $\mathbb{Z}_2 = \langle h \rangle$, where $h$ is the map $\alpha \mapsto \alpha^{-1}$.
\end{abstract}




%


\section{Introduction}

Let $X=(V,E)$ be a simple, undirected graph.  An automorphism of $X$ is a permutation of its vertex set that preserves adjacency (cf. Tutte \cite{Tutte:1966}, Biggs \cite{Biggs:1993}). The set $\{g \in \Sym(V): E^g=E\}$ of all automorphisms of $X$ is called the automorphism group of $X$, and is denoted by $\Aut(X)$.  Given a group $H$ and a subset $S \subseteq H$ such that $1 \notin S$ and $S=S^{-1}$, the Cayley graph of $H$ with respect to the $S$, denoted by $\Cay(H,S)$, is defined to be the graph with vertex set $H$ and edge set $\{(h,sh): h \in H, s \in S\}$. The right regular representation $R(H)$ acts as a group of automorphisms of the Cayley graph  $\Cay(H,S)$, and hence a Cayley graph is always vertex-transitive.  The set of automorphisms of the group $H$ that fixes $S$ setwise, denoted by  $\Aut(H,S)$, is a subgroup of the stabilizer $\Aut(\Cay(H,S))_e$ of the vertex $e$ (cf. \cite{Biggs:1993}).  A Cayley graph $\Cay(H,S)$ is said to be normal if $R(H)$ is a normal subgroup of $\Aut(\Cay(H,S))$, or 
equivalently, if $\Aut(\Cay(H,S)) = R(H) \rtimes \Aut(H,S)$ (cf. \cite{Godsil:1981}, \cite{Xu:1998}). 

An open problem in the literature is to determine which Cayley graphs are normal.  Let $S$ be a set of transpositions generating $S_n$.  The transposition graph of $S$ is defined to be the graph with vertex set $\{1,\ldots,n\}$, and with two vertices $i$ and $j$ being adjacent in this graph iff $(i,j) \in S$.  A set of transpositions $S$ generates $S_n$ iff the transposition of $S$ is connected.  Godsil and Royle \cite{Godsil:Royle:2001} showed that if the transposition graph of $S$ is an asymmetric tree, then $\Cay(S_n,S)$ has automorphism group isomorphic to $S_n$.  Feng \cite{Feng:2006} showed that if the transposition graph of $S$ is an arbitrary tree, then $\Cay(S_n,S)$ has automorphism group $R(S_n) \rtimes \Aut(S_n,S)$.  Ganesan \cite{Ganesan:DM:2013} showed that if the girth of the transposition graph of $S$ is at least 5, then $\Cay(S_n,S)$ has automorphism group $R(S_n) \rtimes \Aut(S_n,S)$.  In all these cases, the Cayley graph $\Cay(S_n,S)$ is normal.  Ganesan \cite{Ganesan:DM:2013} showed that 
if 
the transposition graph of $S$ is a 4-cycle graph, then $\Cay(S_n,S)$ is not normal.

While one can often obtain some automorphisms of a graph, it is often difficult to prove that one has obtained the (full) automorphism group.  In the present paper, we obtain the full automorphism group of the complete transposition graph.  The complete transposition graph has also been studied for consideration as the topology of interconnection networks \cite{Stacho:Vrto:1998}.  Many authors have investigated the automorphism group of other graphs that arise as topologies of interconnection networks; for example, see \cite{Deng:Zhang:2011}, \cite{Deng:Zhang:2012}, \cite{Zhou:2011}, \cite{Zhang:Huang:2005}.

\textbf{Notation.} Throughout this paper, $S$ represents a set of transpositions generating $S_n$, $X:=\Cay(S_n,S)$ and $G:=\Aut(\Cay(S_n,S))$. $X_r(e)$ denotes the set of vertices in $X$ whose distance to the identity vertex $e$ is exactly $r$.  Thus, $X_0(e) = \{e\}$ and $X_1(e)=S$.  Greek letters $\alpha, \beta,\ldots \in S_n$ usually represent the vertices of $X$ and lowercase Latin letters $g, h,\ldots \in \Sym(S_n)$ often represent automorphisms of $X$.  The support of a permutation $\alpha$ is the set of points moved by $\alpha$. For a graph $X$, $L_e:=L_e(X)$ denotes the set of automorphisms of $X$ that fixes the vertex $e$ and each of its neighbors in $X$.

\bigskip The main result of this paper is the following:

\begin{Theorem}
 Let $S$ be the set of all transpositions in $S_n$ ($n \ge 3$).  Then the automorphism group of the complete transposition graph $\Cay(S_n,S)$ is 
  \[ \Aut(\Cay(S_n,S)) = (R(S_n) \rtimes \Inn(S_n)) \rtimes \mathbb{Z}_2, \]
where $R(S_n)$ is the right regular representation of $S_n$, $\Inn(S_n)$ is the inner automorphism group of $S_n$, and $\mathbb{Z}_2 = \langle h \rangle$, where $h$ is the map $\alpha \mapsto \alpha^{-1}$. The complete transposition graph  $\Cay(S_n,S)$ is not normal.
\end{Theorem}

\section{Preliminaries}

Whitney \cite{Whitney:1932} investigated whether a graph $T$ is uniquely determined by its line graph $L(T)$ and showed that the answer is in the affirmative for all connected graphs $T$ on 5 or more vertices (this is because the only exceptions occur when $T$ is $K_3$ or $ K_{1,3}$, which have fewer than 5 vertices). More specifically, two connected graphs on 5 or more vertices are isomorphic iff their line graphs are isomorphic.  And if $T$ is a connected graph that has 5 or more vertices, then every automorphism of the line graph $L(T)$ is induced by a unique automorphism of $T$, and the automorphism groups of $T$ and of $L(T)$ are isomorphic:

\begin{Theorem} \label{thm:Whitney:graph:linegraph:sameautgroup} (Whitney \cite{Whitney:1932}) 
Let $T$ be a connected graph containing at least 5 vertices.  Then the automorphism group of $T$ and of its line graph $L(T)$ are isomorphic.   
\end{Theorem}

\begin{Theorem} \label{thm:Feng:Aut:Sn:S:equals:AutTS} (Feng \cite{Feng:2006}) Let $S$ be a set of transpositions in $S_n$, and let $T=T(S)$ denote the transposition graph of $S$. Then,
$\Aut(S_n,S) \cong \Aut(T)$.
\end{Theorem}

Feng's result (Theorem~\ref{thm:Feng:Aut:Sn:S:equals:AutTS}) does not require that $S$ generate $S_n$, i.e. it holds even if the transposition graph of $S$ is not connected.

\begin{Theorem} \label{thm:Aut:Sn:S:equals:Inn:Sn} (Suzuki \cite[Chapter 3, Section 2]{Suzuki:1982}
If $n \ge 2$ and $n \ne 6$, then $\Aut(S_n)=\Inn(S_n)$.  If $n=6$, then $|\Aut(S_n):\Inn(S_n)|=2$, and every element in $\Aut(S_n) - \Inn(S_n)$ maps a transposition to a product of three disjoint transpositions.
\end{Theorem}


\section{An equivalent condition for normality}

Let $S$ be a set of transpositions generating $S_n$ ($n \ge 5$).  Let $X:=\Cay(S_n,S)$ and let $L_e=L_e(X)$ denote the set of automorphisms of $X$ that fixes the identity vertex $e$ and each of its neighbors.  In this section an equivalent condition for normality of $\Cay(S_n,S)$ is obtained: the Cayley graph $\Cay(S_n,S)$ is normal iff $L_e=1$.  It is not assumed in this section that $S$ is the complete set of transpositions in $S_n$.

\begin{Lemma} \label{lemma:uniqueC4}
Let $S$ be a set of transpositions generating $S_n$.  Let $\tau, \kappa \in S, \tau \ne \kappa$. Then, $\tau \kappa = \kappa \tau$ if and only if there is a unique 4-cycle in $\Cay(S_n,S)$ containing $e,\tau$ and $\kappa$.
\end{Lemma}

\noindent \emph{Proof}: Suppose $\tau \kappa=\kappa\tau$.  Then $\tau$ and $\kappa$ have disjoint support.  Let $\omega$ be a common neighbor of the vertices $\tau$ and $\kappa$ in the Cayley graph $\Cay(S_n,S)$.  By definition of the adjacency relation in the Cayley graph, there exist $x,y\in S$ such that $x\tau=y\kappa=\omega$. Observe that $x\tau=y\kappa$ iff $\tau\kappa=xy$.  But since $\kappa$ and $\tau$ have disjoint support, $\tau\kappa=xy$ iff $\tau=x$ and $\kappa=y$ or $\tau=y$ and $\kappa=x$.  Thus, $\omega$ is either the vertex $e$ or the vertex $\tau\kappa$.  Hence, there exists a unique 4-cycle in $\Cay(S_n,S)$ containing $e,\tau$ and $\kappa$, namely the cycle $(e,\tau,\tau\kappa=\kappa\tau,\kappa,e)$.  

To prove the converse, suppose $\tau\kappa \ne \kappa\tau$.  Then $\tau$ and $\kappa$ have overlapping support; without loss of generality, take $\tau=(1,2)$ and $\kappa=(2,3)$.  We consider two cases, depending on whether $(1,3) \in S$. First suppose $(1,3) \notin S$. Let $\omega$ be a common neighbor of $\tau$ and $\kappa$.  So $\omega=x\tau=y\kappa$ for some $x,y \in S$.  As before, $x\tau=y\kappa$ iff $xy=\tau\kappa=(1,2)(2,3)=(1,3,2)$.  The only ways to decompose $(1,3,2)$ as a product of two transpositions is as $(1,3,2)=(1,2)(2,3)=(3,2)(1,3)=(1,3)(1,2)$.  Since $(1,3) \notin S$, we must have $x=(1,2)$ and $y=(2,3)$, whence $\omega=e$.  Thus, $\tau$ and $\kappa$ have only one common neighbor, namely $e$. Therefore, there does not exist any 4-cycle in $\Cay(S_n,S)$ containing $e,\tau$ and $\kappa$.

Now suppose $\rho:=(1,3) \in S$.  Then $S$ contains the three transpositions $\tau=(1,2),\kappa=(2,3)$ and $\rho=(1,3)$.  The Cayley graph of the permutation group generated by these transpositions is the complete bipartite graph $K_{3,3}$. Hence $\Cay(S_n,S)$ contains as a subgraph the complete bipartite graph $K_{3,3}$ with bipartition $\{e,\kappa\tau,\tau\kappa \}$ and $\{\tau,\kappa,\rho \}$.  There are exactly two  4-cycles in $\Cay(S_n,S)$ containing $e,\tau$ and $\kappa$, namely the 4-cycle through the vertex $\kappa\tau$ and the 4-cycle through the vertex $\kappa\tau$. Thus, while there exists a 4-cycle in this case, it is not unique.
\qed

\begin{Proposition} \label{prop:autHS:in:autLT}
Let $S$ be a set of transpositions generating $S_n$.  Then, every automorphism of $S_n$ that fixes $S$ setwise, when restricted to $S$, is an automorphism of the line graph of the transposition graph of $S$. 
\end{Proposition}

\noindent \emph{Proof}: Let $g \in \Aut(S_n,S)$. Let $\tau,\kappa \in S, \tau \ne \kappa$. Since $g$ is an automorphism of $S_n$, it takes $\tau \kappa$ to $(\tau \kappa)^g = \tau^g \kappa^g$.  An automorphism of a group preserves the order of the elements, whence $\tau$ and $\kappa$ have disjoint support if and only if $\tau^g$ and $\kappa^g$ have disjoint support.  Since $g$ fixes $S$, $\tau^g, \kappa^g \in S$.  Thus, in the transposition graph of $S$, the edges $\tau$ and $\kappa$ are incident to a common vertex if and only if the edges $\tau^g$ and $\kappa^g$ are incident to a common vertex.  In other words, $g$ restricted to $S$ is an automorphism of the line graph of the transposition graph of $S$.\qed

Since $\Aut(S_n,S) \subseteq G_e$, a stronger result than Proposition~\ref{prop:autHS:in:autLT} is the following:

\begin{Proposition}  \label{prop:Ge:restrictedtoS:is:in:AutLT}
Let $S$ be a set of transpositions generating $S_n$, and let $G:=\Aut(\Cay(S_n,S))$.  If $g \in G_e$, then $g$ restricted to $S$ is an automorphism of the line graph of the transposition graph of $S$.
\end{Proposition}

\noindent \emph{Proof}: Let $\tau, \kappa \in S$ and $g \in G_e$.  Let $L(T)$ denote the line graph of the transposition graph of $S$.  Two transpositions commute iff they have disjoint support.  It needs to be shown that the restriction of $g$ to $S$ is an automorphism of $L(T)$, i.e. that $\tau,\kappa$ have disjoint support iff $\tau^g,\kappa^g$ have disjoint support.  Thus, it suffices to show that $\tau \kappa=\kappa\tau$ iff $\tau^g \kappa^g=\kappa^g \tau^g$.  By Lemma~\ref{lemma:uniqueC4}, $\tau \kappa=\kappa\tau$ iff there is a unique 4-cycle in $\Cay(S_n,S)$ containing $e, \tau$ and $\kappa$, which is the case iff there is a unique 4-cycle in $\Cay(S_n,S)$ containing $e, \tau^g$ and $\kappa^g$, which is the case iff $\tau^g \kappa^g = \kappa^g \tau^g$.
\qed

\begin{Proposition} \label{prop:restriction:map:is:surjective}
 Let $S$ be a set of transpositions generating $S_n$, and let $T=T(S)$ denote the transposition graph of $S$ and $L(T)$ denote its line graph.  Let $G:=\Aut(\Cay(S_n,S))$.  Then the restriction map from $G_e$ to $\Aut(L(T))$ defined by $g \mapsto g |_S$ is surjective.
\end{Proposition}

\noindent \emph{Proof}: Let $h \in \Aut(L(T))$. Then $h \in \Sym(S)$ since the vertices of $L(T)$ correspond to the transpositions in $S$.  We show that there exists an element $g \in G_e$ whose action on $S$ is identical to that of $h$. By Whitney's Theorem~\ref{thm:Whitney:graph:linegraph:sameautgroup}, there is an automorphism $h' \in \Aut(T)$ that induces $h$.  Now $h'$ is a permutation in $S_n$. Let $g$ denote conjugation by $h'$.  Thus, $g \in \Aut(S_n)$.  Since $h'$ is an automorphism of $T$, it fixes the edge set $S$ of $T$.  Hence conjugation by $h'$ also fixes $S$, i.e., $g \in \Aut(S_n,S)$.  Since $\Aut(S_n,S) \subseteq G_e, g \in G_e$.  It is clear that $g|_S$ equals $h$.  For example, if $h$ takes $\{i,j\}$ to $\{m,\ell\}$, then there exists an $h' \in S_n$ that takes $\{i,j\}$ to $\{i^{h'},j^{h'} \} = \{m,\ell\}$.  Then $g$ takes $(i,j) \in S$ to $(m,\ell) \in S$.  Thus, $g|_S$ and $h$ induce the same permutation of $S$, which implies the given restriction map is surjective.
\qed

\begin{Theorem} \label{thm:normal:iff:Le:equals:1}
 Let $S$ be a set of transpositions generating $S_n$ ($n \ge 5$).  Let $L_e$ denote the set of automorphisms of the Cayley graph $\Cay(S_n,S)$ that fixes the vertex $e$ and each of its neighbors.  Then, $\Cay(S_n,S)$ is normal if and only if $L_e=1$. 
\end{Theorem}

\noindent \emph{Proof}: $\Leftarrow$:  Consider the map $f$ from the domain $G_e$, defined to be the restriction map $g \mapsto g|_S$.  By Proposition~\ref{prop:Ge:restrictedtoS:is:in:AutLT}, $f$ is into $\Aut(L(T))$.  The kernel of the map $f: G_e \rightarrow \Aut(L(T))$ is the set of elements in $G_e$ that fixes each element in $S$ and hence equals $L_e$.  Since $L_e=1$, $f$ is injective.  By Proposition~\ref{prop:restriction:map:is:surjective}, $f$ is surjective.  The restriction map is also a homomorphism.  Hence $f$ is an isomorphism.

Thus, $|G_e| = |\Aut(L(T))|$.  By Whitney's Theorem~\ref{thm:Whitney:graph:linegraph:sameautgroup}, the transposition graph $T$ and its line graph $L(T)$ have isomorphic automorphism groups.  Thus, $|G_e| = |\Aut(T)|$.  By Theorem~\ref{thm:Feng:Aut:Sn:S:equals:AutTS}, $\Aut(T) \cong \Aut(S_n,S)$.  Thus, $G_e = \Aut(S_n,S)$, which implies $\Cay(S_n,S)$ is normal.

$\Rightarrow$: If $\Cay(S_n,S)$ is normal, then $G_e = \Aut(S_n,S)$ (cf. \cite{Xu:1998}).  Once again, by Theorem~\ref{thm:Whitney:graph:linegraph:sameautgroup} and Theorem~\ref{thm:Feng:Aut:Sn:S:equals:AutTS}, $|G_e| = |\Aut(L(T))|$.  Also, the map $f: G_e \rightarrow \Aut(L(T)), g \mapsto g|_S$ is surjective by Proposition~\ref{prop:restriction:map:is:surjective}.  Thus, $f$ is also injective and therefore its kernel $L_e=1$.
\qed


\section{Non-normality of the complete transposition graph}

\begin{Proposition} \label{prop:inverse:map:is:aut}
Let $S$ be the set of all transpositions in $S_n~ (n \ge 3)$.  Then, the map $\alpha \mapsto \alpha^{-1}$ is an automorphism of the complete transposition graph $\Cay(S_n,S)$.
\end{Proposition}

\noindent \emph{Proof}: 
Let $G$ denote the automorphism group of the Cayley graph $\Cay(S_n,S)$ and let $e$ denote the identity element in $S_n$.  The Cayley graph $\Cay(S_n,S)$ is normal if and only if the stabilizer $G_e \subseteq \Aut(S_n)$ (cf. Xu \cite[Proposition 1.5]{Xu:1998}).  Thus, to prove that $\Cay(S_n,S)$ is not normal, it suffices to show that $G_e$ contains an element which is not a homomorphism from $S_n$ to itself.  Consider the map $\alpha \mapsto \alpha^{-1}$ from $S_n$ to itself.  Since $n \ge 3$, $S_n$ is nonabelian, whence the map $\alpha \mapsto \alpha^{-1}$ is not a homomorphism.  It suffices to show that the map $\alpha \mapsto \alpha^{-1}$ is an automorphism of the Cayley graph $\Cay(S_n,S)$.  

Let $\alpha$ and $\beta$ be two adjacent vertices in the graph $\Cay(S_n,S)$.  Then $\alpha$ and $\beta$ differ by a transposition, i.e. there is some $i \ne j$ such that $\beta=(i,j)g$.  We shall prove that $\alpha^{-1}$ and $\beta^{-1}$ also differ by a transposition; since the set $S$ contains all transpositions in $S_n$, it follows that $\alpha^{-1}$ and $\beta^{-1}$ are also adjacent vertices in $\Cay(S_n,S)$.

Two cases arise, depending on whether $i$ and $j$ are in the same cycle of $\alpha$ or in different cycles of $\alpha$. Suppose $i$ and $j$ are in the same cycle of $\alpha$, say $\alpha=(\alpha_1,\ldots,\alpha_r,i,\beta_1,\ldots,\beta_s,j) \cdots$.  Then $\beta=(i,j)\alpha=(\alpha_1,\ldots,\alpha_r,i,\beta_1,\ldots,\beta_s,j) \cdots$.  A quick calculation shows that $\alpha^{-1}$ and $\beta^{-1}$ differ by the transposition $\tau=(\alpha_1,\beta_1)$ if $r,s \ge 1$, by $\tau=(i,\beta_1)$ if $r=0, s \ge 1$,  by $\tau=(j,\alpha_1)$ if $s=0, r \ge 1$, and by $\tau=(i,j)$ if $r,s=0$.  Hence $\alpha^{-1}=\tau h^{-1}$ for some transposition $\tau$.  Thus, $\alpha^{-1}$ and $\beta^{-1}$ are also adjacent vertices in the Cayley graph $\Cay(S_n,S)$.

Suppose $i$ and $j$ are in different cycles of $\alpha$, say $\alpha=(\alpha_1,\ldots,\alpha_r,i)(\beta_1,\ldots,\beta_s,j) \cdots$ and $\beta=(i,j) \alpha$. Then $i$ and $j$ are in the same cycle of $\beta$ and $(i,j) \beta= \alpha$.  By the argument in the previous paragraph applied to $\beta$ instead of $\alpha$, it follows that $\beta^{-1}$ and $\alpha^{-1}$ are adjacent vertices in $\Cay(S_n,S)$.

We have shown that if $\alpha$ and $\beta$ are adjacent vertices, then so are $\alpha^{-1}$ and $\beta^{-1}$.  It follows that if $\alpha^{-1}$ and $\beta^{-1}$ are adjacent vertices, then so are $(\alpha^{-1})^{-1}=\alpha$ and $(\beta^{-1})^{-1}=\beta$.  Hence, $\alpha \mapsto \alpha^{-1}$ is an automorphism of the Cayley graph $\Cay(S_n,S)$.
\qed

\begin{Theorem} \label{thm:aut:completetransp:subgroup}
 Let $S$ be the set of all transpositions in $S_n$ ($n \ge 3$).  Then 
 \[ \Aut(\Cay(S_n,S)) \supseteq (R(S_n) \rtimes \Inn(S_n)) \rtimes \mathbb{Z}_2, \]
 where $R(S_n)$ is the right regular representation of $S_n$, $\Inn(S_n)$ is the inner automorphism group of $S_n$, and $\mathbb{Z}_2 = \langle h \rangle$, and $h$ is the map $\alpha \mapsto \alpha^{-1}$.
\end{Theorem}

\noindent \emph{Proof}: 
Let $G:=\Aut(\Cay(S_n,S))$ denote the automorphism group of the complete transposition graph.  Since $R(S_n)$ and $\Aut(S_n,S)$ are automorphisms of the Cayley graph (cf. \cite{Biggs:1993}), we have $G \supseteq R(S_n) \rtimes \Aut(S_n,S)$.  Also, $S$ is a nonempty set of transpositions, so by Theorem~\ref{thm:Aut:Sn:S:equals:Inn:Sn} every element in $\Aut(S_n,S)$ is an inner automorphism of $S_n$.  In fact, the elements in $\Aut(S_n,S)$ are exactly conjugations by the automorphisms of the transposition graph of $S$ (cf. Theorem~\ref{thm:Feng:Aut:Sn:S:equals:AutTS} and \cite{Feng:2006}).  The transposition graph of $S$ is complete, hence $\Aut(S_n,S) = \Inn(S_n) \cong S_n$.  

By Proposition~\ref{prop:inverse:map:is:aut} the map $(h: \alpha \mapsto \alpha^{-1})$ is in $G$.  We show that  $h \notin R(S_n) \rtimes \Inn(S_n)$.  By way of contradiction, suppose $h \in R(S_n) \rtimes \Inn(S_n)$.  Then $h=ab$ for some $a \in R(S_n), b \in \Inn(S_n)$.  Hence $e^h=e^{-1}=e$, and $e^{ab}=e$.  Since $b \in \Inn(S_n)$, $b$ fixes $e$.  Thus $e^a=a$, whence $a=1$.  Thus, $h=ab=b \in \Inn(S_n)$, which is a contradiction since the map $h: \alpha \mapsto \alpha^{-1}$ is not a homomorphism.   

Thus $G$ contains $H:= (R(S_n) \rtimes \Inn(S_n)) \rtimes \mathbb{Z}_2$, where $\mathbb{Z}_2 := \langle h \rangle$ and $R(S_n) \rtimes \Inn(S_n)$ has index 2 in $H$ and hence is a normal subgroup in $H$.
\qed

This implies that the complete transposition graph $\Cay(S_n,S)$ has at least $2(n!)^2$ automorphisms, for all $n \ge 3$. 

\begin{Theorem}
 Let $S$ be the set of all transpositions in $S_n$ ($n \ge 3$).  Then the complete transposition graph $\Cay(S_n,S)$ is not normal.
\end{Theorem}

\noindent \emph{First proof}: By Proposition~\ref{prop:inverse:map:is:aut}, the inverse map $h: \alpha \mapsto \alpha^{-1}$ is an automorphism of the Cayley graph $\Cay(S_n,S)$.  The map $h$ fixes the vertex $e$ and also fixes each transposition $(i,j) \in S$.  Thus, $h \in L_e$.  Since $n \ge 3$, $\exists \alpha \in S_n$ such that $\alpha \ne \alpha^{-1}$.  Thus $h$ is not the trivial map and $L_e > 1$.  If $n=3$ or $n=4$, it can be confirmed through computer simulations that $R(S_n)$ is not a normal subgroup of the automorphism group of $\Cay(S_n,S)$; hence $\Cay(S_n,S)$ is not normal in these cases.  If $n \ge 5$, then Theorem~\ref{thm:normal:iff:Le:equals:1} applies and again $\Cay(S_n,S)$ is not normal.\qed

\bigskip \noindent \emph{Second proof}: Alternatively, Theorem~\ref{thm:aut:completetransp:subgroup} provides a second proof that the complete transposition graph  is not normal:  a normal Cayley graph $\Cay(S_n,S)$ has the smallest possible full automorphism group $R(S_n) \rtimes \Aut(S_n,S)$, whereas by Theorem~\ref{thm:aut:completetransp:subgroup} the complete transposition graph has an automorphism group that is strictly larger.  Hence the complete transposition graph is not normal.
\qed

Let $S$ be a set of transpositions generating $S_n$ ($n \ge 3$).  The only Cayley graphs $\Cay(S_n,S)$ known so far to be non-normal are those arising from the 4-cycle transposition graph and from the transposition graphs that are complete.

\section{Automorphism group of the complete transposition graph}

Let $S$ be the set of all transpositions in $S_n$.  In the previous section a set of $2(n!)^2$ automorphisms were exhibited for the complete transposition graph $\Cay(S_n,S)$.  In this section, it is proved that the complete transposition graph has no other automorphisms, which implies that the subgroup given in Theorem~\ref{thm:aut:completetransp:subgroup} is in fact the full automorphism group.

\begin{Theorem} \label{thm:Le:equals:C2}
 Let $S$ be the set of all transpositions in $S_n$ and let $X$ be the complete transposition graph $\Cay(S_n,S)$.  Let $L_e(X)$ denote the set of automorphisms of $X$ that fixes the vertex $e$ and each of its neighbors.  Then $L_e(X) = \{1,h \}$, where $h: V(X) \rightarrow V(X)$ is the map $\alpha \rightarrow \alpha^{-1}$.
\end{Theorem}

\noindent \emph{Proof}:
By Proposition~\ref{prop:inverse:map:is:aut}, $L_e \supseteq \{1,h\}$.  We need to show that $L_e$ has no other elements. 

The vertex $e$ of $X$ corresponding to the identity element in $S_n$ has as its neighbors the set $S$ of all  transpositions in $S_n$.  Suppose $g$ is an automorphism of $X$ that fixes the vertex $e$ and each vertex in $S$; so $g \in L_e(X)$.  Then the set of common neighbors of the three vertices $(1,2), (2,3)$ and $(1,3)$ in $S$, namely the set $\Delta:=\{(1,2,3),(1,3,2)\}$, is a fixed block of $g$.  We show that the action of $L_e:=L_e(X)$ on $\Delta$ uniquely determines its action on all the remaining vertices, i.e. that if $g \in L_e$ fixes $\Delta$ pointwise, then $g=1$, and if $g$ interchanges $(1,2,3)$ and $(1,3,2)$, then $g$ extends uniquely to the automorphism $\alpha \mapsto \alpha^{-1}$ of $X$.

\begin{figure}
\centering
\begin{tikzpicture}[scale=1]
\vertex[fill] (e) at (0,0) [label=below:$e$] {};

\vertex[fill] (v12) at (2,3) [label=above:$(12)$] {};
\vertex[fill] (v23) at (2,2) [label=above:$(23)$] {};
\vertex[fill] (v13) at (2,1) [label=above:$(13)$] {};

\vertex[fill] (alpha) at (4,4) [label=above:{$\alpha=(123)$}] {};
\vertex[fill] (alphainv) at (4,3) [label=below:{$\alpha^{-1}$}] {};
\vertex[fill] (beta) at (4,2) [label=below:{$\beta$}] {};
\vertex[fill] (betainv) at (4,1) [label=below:{$\beta^{-1}$}] {};

\draw (6,3.5) ellipse (.2cm and 1cm);
\node at (6,2.2) {$W_\alpha$};
\node (walphatop) at (6.1,4.5) {};
\node (walphabottom) at (6.1,2.4) {};
\path (alpha) edge (walphatop);
\path (alpha) edge (walphabottom);

\node at (2,-1) {$\underbrace{~~~~~~}_{X_1(e)}$};
\node at (4,-1) {$\underbrace{~~~~~~}_{X_2(e)}$};
\node at (6,-1) {$\underbrace{~~~~~~}_{X_3(e)}$};
\path
    (e) edge (v12)
    (e) edge (v13)
    (e) edge (v23)
    (v12) edge (alpha)
    (v23) edge (alpha)
    (v13) edge (alpha)
    (v12) edge (alphainv)
    (v23) edge (alphainv)
    (v13) edge (alphainv)
;
\end{tikzpicture}
\caption{Distance partition of the Cayley graph $\Cay(S_n,S)$} 
\label{fig:distance:partition}
\end{figure}
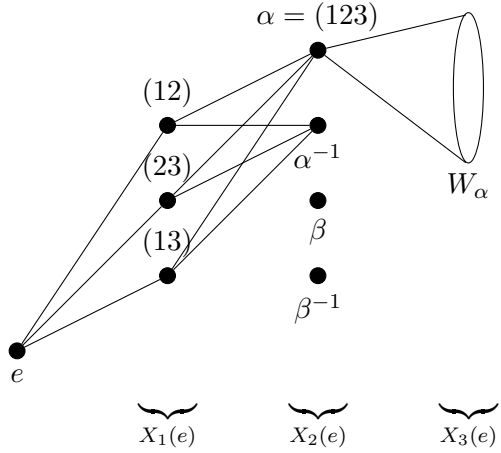

Suppose $g \in L_e$ and $g$ fixes $\Delta = \{\alpha,\alpha^{-1} \}$ pointwise, where $\alpha=(1,2,3)$.  Let $\beta=(2,3,4)$.  We show $g$ fixes $\{\beta,\beta^{-1} \}$ also pointwise.  Given any vertex $\gamma \in V(X)$ that is a 3-cycle permutation (so the distance in $X$ between $\gamma$ and $e$ is 2), let $W_\gamma$ be the set of neighbors of $\gamma$ that have distance 3 to $e$ in $X$ (see Figure~\ref{fig:distance:partition}). 

We claim that $|W_\alpha \cap W_\beta|=|W_{\alpha^{-1}} \cap W_{\beta^{-1}}|=2$ and $|W_\alpha \cap W_{\beta^{-1}} | = |W_{\alpha^{-1}} \cap W_\beta|=1$.  Supose some neighbor of $\alpha=(1,2,3)$ is also a neighbor of $\beta=(2,3,4)$.  Then $\exists x, y \in S$ such that $x \alpha = y \beta$.  Hence, $\alpha \beta^{-1} = (1,2,3)(2,3,4)=(1,4,3)=x^{-1}y=xy$.  Now $(1,4,3)=(1,4)(1,3)=(1,3)(3,4)=(3,4)(1,4)$.  So $x \in \{(1,4),(1,3),(3,4)\}$.  But if $x=(1,3)$, then $x \alpha = (1,3)(1,2,3)=(2,3)$, so $x \alpha$ has distance 1 to $e$, and $x \alpha \notin W_\alpha$.  Thus, there are two solutions $(1,4)$ and $(3,4)$ for $x$ in $x \alpha=y \beta \in W_\alpha \cap W_\beta$.  Hence $|W_\alpha \cap W_\beta|=2$.  Similarly, $|W_{\alpha^{-1}} \cap W_{\beta^{-1}}|=2$.  Now consider  $|W_\alpha \cap W_{\beta^{-1}} |$.  If $x,y \in S$ are such that $x \alpha=y \beta^{-1}$, then $\alpha \beta=(1,2,3)(2,3,4)=(1,3)(2,4)=xy$.  But if $x=(1,3)$, then $y \alpha=(1,3)(1,2,3)=(2,3) \notin W_\alpha$.  Thus, $x=(2,4)$, $y=(1,3)$, 
$x \alpha = (2,4)(1,2,3)=(1,2,4,3)$, and  $|W_\alpha \cap W_{\beta^{-1}} | = |\{(1,2,4,3)\}|=1$.

Since $g$ is an automorphism of $X$, it preserves the number of common neighbors of any two vertices.  Thus, if $g$ fixes $\alpha$ and $\alpha^{-1}$, by the result in the previous paragraph, $g$ also fixes $\beta$ and $\beta^{-1}$.  More generally, if $g$ fixes vertex $(j,k,i)$, then $g$ also fixes $(j,k,\ell)$ for each $\ell \neq j,k,i$.  Repeating this process, we see that $g$ fixes all vertices that are 3-cycles in $S_n$.  The only other vertices having distance 2 to $e$ in $X$ are those permutations that are a product of two disjoint transpositions, and each of these vertices are also fixed by $g$ by Lemma~\ref{lemma:uniqueC4}. 

Thus, if $g \in L_e(X)$ fixes vertex $(1,2,3)$, then $g$ fixes each vertex that has distance 2 to $e$.  Let $X_r(e)$ denote the set of vertices that have distance $r$ to $e$. We have that $g$ fixes $X_0(e)$ and $X_1(e)$ pointwise since $g \in L_e$, and it was just shown that if $g$ fixes $(1,2,3) \in X_2(e)$, then $g$ also fixes $X_2(e)$ pointwise.  Since $g$ is an automorphism, it maps the neighbors of a vertex $\alpha$ to the neighbors of $\alpha^g$.  But by the next proposition (Proposition~\ref{prop:distinct:neighbors}), any two distinct vertices in $X_k(e)$ $(k \ge 3)$ have a different set of neighbors in $X_{k-1}(e)$.  Thus, if $g$ fixes $X_{k-1}(e)$ pointwise, then $g$ also fixes $X_k(e)$ pointwise.  By induction on $k$, $g$ is the trivial automorphism.

If $g \in L_e$ interchanges $(1,2,3)$ and $(1,3,2)$, and $h$ is the map $\alpha \mapsto \alpha^{-1}$, then $gh=1$ by the previous paragraph, whence $g=h^{-1}=h$.  Thus, $L_e = \{1,h\} \cong C_2$.
\qed

\bigskip In the proof above, we used the following result:

\begin{Proposition} \label{prop:distinct:neighbors}
 Let $n \ge 5$ and let $X=\Cay(S_n,S)$ be the complete transposition graph.  Let $\alpha$ and $\beta$ be distinct vertices in $X_k(e)$ $(k \ge 3)$.  Then the set of neighbors of $\alpha$ in $X_{k-1}(e)$ and of $\beta$ in $X_{k-1}(e)$ are not equal. 
\end{Proposition}

\noindent \emph{Proof}:
Each permutation in $X_k(e)$ can be written as a product of $k$ transpositions, and since the length of this product is minimal, the edges of the transposition graph of $S$ corresponding to these $k$ transpositions form a forest.

Let $\alpha, \beta \in X_k(e)$. If the support of $\alpha$ and of $\beta$ are not equal, then they clearly have different sets of neighbors in $X_{k-1}(e)$ because some transposition in a forest that yielded $\alpha$ is incident to a vertex that does not belong to any forest that yields $\beta$. (For example, if $\alpha=(1,2,3)(4,5)$ and $\beta=(1,2,3,4)$ are two vertices in $X_3(e)$, then $\alpha$ does have a neighbor $(1,2)(4,5)$ in $X_2(e)$ whose support contains 5, but $\beta$ does not have such a neighbor.)

Now suppose $\alpha$ and $\beta$ are distinct vertices in $X_k(e)$ that have the same support.  Since $\alpha \ne \beta$, there is a point in their common support, 1 say, such that $1^\alpha \ne 1^\beta$.  So suppose $\alpha=(1,2,x_1,\ldots,x_r) \alpha'$ and $\beta=(1,3,y_1,\ldots,y_t) \beta'$.  We consider three cases:

Case 1:  Suppose $\alpha' = \beta'=1$.  Then $\alpha$ and $\beta$ are cyclic permutations of the same length $r$, where $r \ge 4$ since $k \ge 3$.  If $\alpha=\beta^{-1}$, then we can find two consecutive points in the cycle of $\alpha$ that are not consecutive in the cycle of $\beta$.  Suppose $i,j$ are these two points; so $\alpha=(i,j,k,\ldots,m)$ and $\beta=(i,\ell,\ldots,j,p,\ldots)$.  Then $\gamma=(i,j)(k,\ldots,m)$ is a neighbor of $\alpha$ in $X_{k-1}(e)$ but not of $\beta$.  For if $s \gamma=\beta$ for some transposition $s$, then $s=\beta \gamma^{-1} = (i,\ell,\ldots,j,p,\ldots)(i,j)(k,m,\ldots)$.  Now $s$ moves $i$ since $i^s = i^{\beta \gamma^{-1}}=\ell^{\gamma^{-1}} \ne i$.  Also, $s$ moves $j$ since $j^s=p^{\gamma^{-1}} \ne j$.  If $s=(i,j)$, then $(i,j) = (i,\ell,\ldots,j,p,\ldots)(k,m,\ldots)(i,j)$, whence $(i,\ell,\ldots,j,p,\ldots,q)(k,m,\ldots)=1$, which is a contradiction since $q$ is not fixed by the left hand side but is fixed by the right hand side.  Thus $s$ moves at least 3 points.  
But 
then $s$ is not a transposition, a contradiction.  Hence $\beta$ does not have $\gamma$ as a neighbor.

If $\alpha=\beta^{-1} = (\alpha_1,\ldots,\alpha_r)$, then $(\alpha_1,\ldots,\alpha_{r-1})$ is a neighbor of $\alpha$ in $X_{k-1}(e)$ but not of $\beta$.  

Case 2:  Suppose $\alpha'=\beta' \ne 1$.  So $\alpha=(1,2,x_1,\ldots,x_r) (\alpha_1,\ldots,\alpha_s) \alpha''$, 
\\ $\beta=(1,3,y_1,\ldots,y_t)(\alpha_1,\ldots,\alpha_s) \alpha''$ for some $s \ge 2$ and some (possibly trivial) permutation $\alpha''$. Let $\gamma=(1,2,x_1,\ldots,x_r)(\alpha_1,\ldots,\alpha_{s-1}) \alpha''$.  Then $\gamma$ is a neighbor of $\alpha$ but not of $\beta$.

Case 3:  Suppose $\alpha' \ne \beta'$.  So $\alpha=(1,2,x_1,\ldots,x_r) \alpha'$ and $\beta = (1,3,y_1,\ldots,y_t) \beta'$ are in $X_k(e)$.  If the support of $\alpha'$ and of $\beta'$ are equal, then take $\gamma:=(1,2,x_1,\ldots,x_r) \gamma'$, where $\gamma'$ is any vertex that is adjacent in $X$ to $\alpha'$ and that lies on a shortest $e-\alpha'$ path in $X$.  Then $\gamma$ is adjacent to $\alpha$ but not to $\beta$.  

On the other hand, if the support of $\alpha'$ and of $\beta'$ are not equal, we consider three subcases:
\\ (i) Suppose $r=t=0$.  Then $\alpha=(1,2)\alpha',\beta=(1,3)\beta'$.  Take $\gamma=(1,2)\gamma'$ where $\gamma'$ is any vertex in $X$ adjacent to $\alpha'$ and such that $\gamma$ lies on a shortest $e-\alpha'$ path in $X$.  Then $\gamma$ is a neighbor of $\alpha$ but not of $\beta$ because if $s$ is a transposition, then $s \gamma = s (1,2) \gamma'$ will either split a cycle in $\gamma$ or merge two cycles in $\gamma$, neither of which can produce $(1,3)\beta'$.
\\ (ii) Suppose $r \ge 1$ and $t=0$.  Then $\beta=(1,3)\beta'$.  Take $\gamma$ to be $(1,2)(x_1,\ldots,x_r) \alpha'$.  As in subcase (i), there does not exist any transposition $s$ such that $s \gamma = \beta$. 
\\ (iii)  Suppose $r,t \ge 1$.  Let $\alpha = (1,2,x_1,\ldots,x_r) \alpha' = \alpha^0 \alpha'$ and $\beta = (1,3,y_1,\ldots,y_t) \beta' = \beta^0 \beta'$.   Let $\supp(\alpha)$ denote the support of the permutation $\alpha$.

If $3 \notin \supp(\alpha^0)$, take $\gamma = (1,3)(y_1,\ldots,y_t) \beta'$.  Then $\gamma$ is a neighbor of $\beta$.  But if $\alpha=s \gamma$ for some transposition $s$, then $s$ must modify the cycle $(1,3)$ of $\gamma$, hence must merge this cycle with another one.  The merged cycle will contain both 1 and 3, whence $s \gamma \ne \alpha$ because $3 \notin \supp(\alpha^0)$. Similarly, if $2 \notin \supp(\beta^0)$, then take $\gamma=(1,2)(x_1,\ldots,x_r) \alpha'$, and $\gamma$ is a neighbor of $\alpha$ but not of $\beta$.  

Finally, suppose $e \in \supp(\alpha^0)$ and $2 \in \supp(\beta^0)$.  Split $\alpha^0=(1,2,\ldots,3,\ldots)$ before the 3 to get $\gamma^0=(1,2,\ldots)(3,\ldots)$.  Let $\gamma:=\gamma^0 \alpha'$.  Then $\gamma$ is a neighbor of $\alpha$.  If $\gamma$ is also a neighbor of $\beta=(1,3,y_1,\ldots,y_t) \beta'$, then $s \gamma=\beta$ for some $s$ that merges the two cycles in $\gamma^0$.  But such a merge will produce a single cycle that has the same support as $\alpha^0$, whereas $\supp(\alpha^0) \ne \supp(\beta^0)$ by hypothesis.  Hence $\gamma$ is not a neighbor of $\beta$.
\qed

\begin{Corollary} \label{cor:ubound:numauts:completetranspgraph}
 Let $S$ be the set of all transpositions in $S_n$ ($n \ge 3$).  Then
 \[ |\Aut(\Cay(S_n,S))| \le 2(n!)^2. \]
\end{Corollary}

\noindent \emph{Proof}:
Let $G:=\Aut(\Cay(S_n,S))$. The upper bound is verified to be exact if $n=3,4$ by computer simulations. If $n \ge 5$, by Lemma~\ref{prop:Ge:restrictedtoS:is:in:AutLT}, every element in $G_e$, when restricted to $S$, is an automorphism of the line graph of the transposition graph of $S$.  The transposition graph is complete, and hence its line graph has automorphism group isomorphic to $S_n$ (cf. Theorem~\ref{thm:Whitney:graph:linegraph:sameautgroup}).  Hence $|G_e| \le |S_n|~ |L_e|$.  Also, $|L_e|=2$, hence $|G_e| \le 2(n!)$.  Thus $|G| = |V(X)|~ |G_e| \le n! (2n!)$.
\qed

By Corollary~\ref{cor:ubound:numauts:completetranspgraph}, subgroup given in Theorem~\ref{thm:aut:completetransp:subgroup} is in fact the full automorphism group:

\begin{Corollary}  Let $S$ be the set of all transpositions in $S_n$ ($n \ge 3$).  Then, the automorphism group of the complete transposition graph $\Cay(S_n,S)$ is 
  \[ \Aut(\Cay(S_n,S)) = (R(S_n) \rtimes \Inn(S_n)) \rtimes \mathbb{Z}_2, \]
where $R(S_n)$ is the right regular representation of $S_n$, $\Inn(S_n)$ is the inner automorphism group of $S_n$, and $\mathbb{Z}_2 = \langle h \rangle$, where $h$ is the map $\alpha \mapsto \alpha^{-1}$.
\end{Corollary}


\bibliographystyle{plain}
\bibliography{refsaut}


\end{document}